\documentclass[a4paper, 11pt]{article}
\usepackage{amsfonts}
\usepackage{amsfonts}
\usepackage{amsfonts}
\usepackage{amsfonts}
\usepackage{mathrsfs}
\usepackage{mathrsfs}
\usepackage{amsfonts,latexsym,rawfonts,amsmath,amssymb,amsthm}

\usepackage{color}
\usepackage[applemac]{inputenc}
\usepackage[T1]{fontenc}
\usepackage{textcomp}
\usepackage{geometry}
\usepackage{graphicx}
\usepackage{mathptmx}
\usepackage{hyperref}

\numberwithin{equation}{section}

\newcommand{\beq}{\begin{equation}}
\newcommand{\eeq}{\end{equation}}
\newcommand{\beqs}{\begin{eqnarray*}}
\newcommand{\eeqs}{\end{eqnarray*}}
\newcommand{\beqn}{\begin{eqnarray}}
\newcommand{\eeqn}{\end{eqnarray}}
\newcommand{\beqa}{\begin{array}}
\newcommand{\eeqa}{\end{array}}

\newtheorem{prop}{Proposition}[section]
\newtheorem{theo}[prop]{Theorem}
\newtheorem{lem}[prop]{Lemma}

\newtheorem{conj}[prop]{Conjecture}

\title{Minimal two-spheres with constant curvature in the quaternionic projective space}
\author{Jie Fei \quad  Chiakuei Peng \quad Xiaowei Xu\footnote{Corresponding author}}

\begin{document}
\bibliographystyle{plain}

\date{}

\maketitle

\noindent\textbf{Abstract}: In this paper we completely classify the homogeneous two-spheres, especially, the minimal homogeneous ones
in the quaternionic projective space $\textbf{HP}^n$. According to our classification, more minimal constant curved two-spheres in $\textbf{HP}^n$ are obtained than Ohnita conjectured in \cite{[Oh]}.

\vspace{0.2cm}

\vspace{0.2cm}

\tableofcontents

\section{Introduction}
The study of minimal constant curved two-spheres in symmetric spaces has a long history. E.Calabi \cite{[Ca]} and do Carmo and Wallach \cite{[DW]} proved that
any isometric minimal immersion from $S^2_K$ into $S^n(1)$ is
congruent to a linearly full one in $S^{2m}(1)$ with
$K=2/m(m\!+\!1)$ for some positive integer $m$. S. Bando, Y.
Ohnita \cite{[BO]} and J. Bolton, G.R.
Jensen, M. Rigoli, L. M. Woodward \cite{[BJRW]} proved that the linearly full minimal $S^2_K$ in
$\textbf{CP}^n$ must be one of the Veronese surfaces with
$K=4/(n\!+\!2j(n\!-\!j))$ for some integer $0\leq\! j\!\leq n$. These perfect works lead to the study of minimal two-spheres with constant curvature in symmetric spaces. For example,
Zh.Q.Li and Zh.H.Yu \cite{[LY]} classified the minimal constant curved two-spheres in the complex Grassmann manifold $G(2,4)$, X.X.Jiao and C.K.Peng \cite{[JP]} classified the holomorphic ones in $G(2,5)$. L.Delisle, V.Hussin and W.J.Zakrzewski proved \cite{[DHZ]} \cite{[DHZ2]} these two results again from the viewpoint of the Grassmannian sigma-model, and they also proposed conjectures w.r.t. the holomorphic constant curved two-spheres in complex Grassmannians. Recently, C.K.Peng and X.W.Xu \cite{[PX]} completely classified the homogenous minimal two-spheres in the complex Grassmann manifold $G(2,n)$, and C.K.Peng, J.Wang and X.W.Xu \cite{[PJX]} completely classified the homogenous minimal ones in the complex hyperquadric $Q_n$.

Naturally, one can consider the corresponding problem when the target space is quaternionic projective space $\textbf{HP}^n$. We call a minimal immersion $f$ from $S^2$ into $\textbf{HP}^n$ is \emph{proper} means that $f(S^2)$ is not contained in any totally geodesic submanifold of $\textbf{HP}^n$. Without loss of generality, we only consider the proper ones in this paper. Base on the works of homogeneous harmonic maps into projective space, Y.Ohnita \cite{[Oh]} proposed the following conjecture:
\begin{conj}
Let $f$ be a proper minimal constant curved immersion from $S^2$ into $\textbf{\emph{HP}}^n$, then $f$ is congruent to one of $\{f_\lambda\;|\;\lambda=1,3,\ldots,2n+1\}$, where $f_\lambda$ is defined in the remark of Proposition $4.1$.
\end{conj}

By using the theory of harmonic sequences, L.He and X.X.Jiao \cite{[HJ]} proved that this conjecture is true for $n=2$ under the totally unramified condition. So far, there is no breakthrough towards this conjecture for general $n>2$. Notice that the homogeneous submanifold is a class of important submanifolds. So, the classification of homogenous two-spheres in quaternionic projective space are helpful in understanding this conjecture. In this paper we completely classify the homogeneous two-spheres in $\textbf{HP}^n$. In particular, we obtain all minimal homogeneous two-spheres in $\textbf{HP}^n$. Namely, we have
\begin{theo}
Let $f:S^2\longrightarrow\textbf{\emph{HP}}^n$ be a proper homogeneous minimal immersion. Then, in terms of homogeneous coordinates, $f$ is congruent to one of the following\emph{:}

\emph{(1)} $f_\lambda:[a,b]\mapsto [\phi_{\lambda,n+1}],$ $\lambda\in\{3,5,\ldots,2n+1\}$, and the Gauss curvature $K=8/[(2n+2)^2-(\lambda^2+1)];$

\emph{(2)} $f_1:[a,b]\mapsto[\phi_{1,n+1}]$ and the Gauss curvature $K=4/[n(n+2)];$

\emph{(3)} $f_{\lambda,m,t}:[a,b]\mapsto[\cos t\phi_{\lambda,m},\texttt{\emph{i}}\sin t\phi_{\lambda,m}]$ for some positive weight $\lambda\in\{3,5,\ldots,n\}$, $t\in(0,\pi/2)$, $2m=n+1$, and the Gauss curvature $K=8/[(n+1)^2-(\lambda^2+1)];$

\emph{(4)} $f_{m_1,m_2}:[a,b]\mapsto[\sqrt{m_1/(m_1+m_2)}\;\phi_{1,m_1},
\texttt{\emph{i}}\sqrt{m_2/(m_1+m_2)}\;\phi_{1,m_2}]$ for some positive $m_1\leq m_2$ so that $m_1+m_2=n+1$ is an even, and the Gauss curvature $K=4/(m_1^2+m_2^2-1);$

\emph{(5)} $f'_{m_1,m_2}:[a,b]\mapsto[\sqrt{m_1/(m_1+m_2)}\;\phi_{1,m_1},
\sqrt{m_2/(m_1+m_2)}\;\phi_{1,m_2}]$ for some positive $m_1<m_2$ so that $m_1+m_2=n+1$ is an odd, and the Gauss curvature $K=4/(m_1^2+m_2^2-1)$,
where $\phi_{\lambda,m}$ is defined in the end of Section 2.1.
\end{theo}
\noindent All these minimal two-spheres obtained are not congruent to each other, up to a rigidity of $\textbf{HP}^n$. So, by comparing Conjecture 1.1 with Theorem 1.2, we obtain three families of minimal two-spheres more than Y.Ohnita conjectured. Particularly, there is a family of minimal two-spheres in $\textbf{HP}^n$ which depend on a parameter $t\in(0,\pi/2)$. This phenomenon appeared in Theorem 1.2 is quite different from the case of complex projective space (see \cite{[BO]}, \cite{[BJRW]}), due to the non-commutativity of quaternion.

In this paper we first show that there exists a quaternionic representation $\rho$ of $SU(2)$ associated to each homogeneous immersion $f$ from $S^2$ into $\textbf{HP}^n$ so that $f(S^2)$ is a $\rho(SU(2))$-orbit in $\textbf{HP}^n$. Then, in order to characterize the geometry of such orbit, we find the "best" base point in it. That is, we prove that such an orbit must contain a point $P_0$
spanned by vectors belonging to a weight space of the quaternionic representation $\rho$. To study the minimal orbits, we obtain a minimality criterion by using the method of moving frame which is inspired from S.S.Chern and J.G.Wolfson's work \cite{[CW2]}.
By using this criterion, up to a rigidity of $\textbf{HP}^n$, we show that the representation
$\rho$ is a direct sum of two different irreducible representations at the most. Thus, we completely classify the homogenous and minimal homogeneous two-spheres in $\textbf{HP}^n$, see Theorem 3.3, Theorem 4.4 respectively.

\section{Preliminaries}
\subsection{The quaternionic representation of $SU(2)$}\label{sec2}
For completeness, we review some facts on the unitary and quaternionic representations of the special unitary group $SU(2)$.

The special unitary group $SU(2)$ is defined by
\begin{equation*}
SU(2)=\left\{g=\left(
          \begin{array}{cc}
            a & b \\
            -\bar{b} & \bar{a} \\
          \end{array}
        \right)
 : |a|^2+|b|^2=1,\;\; a,b\in \textbf{C} \right\}.
\end{equation*}
Its Lie algebra $\mathfrak{su}(2)$ is real spanned by
\begin{equation*}
    \varepsilon_1=\left(
                     \begin{array}{cc}
                        \texttt{i}&0\\
                         0&-\texttt{i}\\
                     \end{array}
                  \right),\;
    \varepsilon_2=\left(
                     \begin{array}{cc}
                         0&1\\
                         -1&0\\
                     \end{array}
                  \right),\;
    \varepsilon_3=\left(
                     \begin{array}{cc}
                         0&\texttt{i}\\
                         \texttt{i}&0\\
                     \end{array}
                  \right),
\end{equation*}
where $\texttt{i}^2=-1$.
The Maurer-Cartan forms of $SU(2)$ is determined by
\begin{equation}\label{2.1}
\Theta:=dgg^{-1}=\left(
                         \begin{array}{cc}
                           \texttt{i}\omega & \varphi\\
                           -\bar{\varphi}  &-\texttt{i}\omega \\
                         \end{array}
                       \right),
\end{equation}
where $\omega$, $\varphi$ are real and complex one-forms
respectively. The Maurer-Cartan equation is given by
\begin{equation*}
d\Theta=\Theta\wedge\Theta,
\end{equation*}
which implies
\begin{equation}\label{2.2}
d(\texttt{i}\omega)=-\varphi\wedge\bar{\varphi},\;\;\;d\varphi=(2\texttt{i}\omega)\wedge\varphi.
\end{equation}
Let
$T=\big\{\mbox{diag}\{e^{\texttt{i}t},e^{\texttt{-i}t}\}\;|\;t\in\mathbb{R}\big\}$
be a maximal torus subgroup of $SU(2)$. It is known that the set
$SU(2)/T=\{[g]=Tg\;|\;g\in SU(2)\}$ is the complex projective space
$\textbf{CP}^1$, which is diffeomorphic to $S^2$. The canonical
metric $\varphi\bar\varphi$ on $\textbf{CP}^1$ has constant curvature $4$ by
(\ref{2.2}).

Let $V_n$ be the $(n\!+\!1)$-dimensional complex vector space of all
complex homogeneous polynomials of degree $n$ w.r.t. the two complex
variables $z_0$ and $z_1$. We define a Hermitian inner product
$(\;,\;)$ on $V_n$ by
\begin{equation*}
( f,g):=\sum\limits_{k=0}^n\;k!(n-k)!\;a_k\bar{b}_k,
\end{equation*}
for $f=\sum\limits_{k=0}^na_kz_0^kz_1^{n-k}$, $g=\sum\limits_{k=0}^nb_kz_0^kz_1^{n-k}\in V_n$. So, we know
$\big\{\textbf{v}_{k,n}=z_0^kz_1^{n-k}/\sqrt{k!(n-k)!}\;\big|\;0\leq k\leq
n\big\}$
is a unitary basis for $V_n$. A unitary representation $\rho_n$ of
$SU(2)$ on $V_n$ is defined by
\begin{equation*}
\rho_n(g)f(z_0,z_1):=f((z_0,z_1)g^{-1})=
f(\bar{a}z_0+\bar bz_1,-bz_0+az_1)
\end{equation*}
for $g\in SU(2)$ and $ f\in V_n$. Under the basis $\{\textbf{v}_{k,n}:0\leq k\leq n\}$, we get a matrix representation $\rho_n:SU(2)\longrightarrow U(n+1)$, $g\mapsto \rho_n(g)$, and $\rho_n(g)$ is described by
$\textbf{v}_{k,n}\rho_n(g):=\sum\limits_{k=0}^n\Lambda_{l\;k}(g)\textbf{v}_{l,n}$,
where
$$\Lambda_{l\;k}(g)=\sqrt{\frac{l!(n-l)!}{k!(n-k)!}}\sum\limits_{p+q=l}
\left(
  \begin{array}{c}
    k\\
    p \\
  \end{array}
\right)\left(
         \begin{array}{c}
           n-k \\
           q \\
         \end{array}
       \right)a^{n-k-q}(-b)^{q}\bar{a}^p\bar{b}^{k-p}
$$
which satisfies $\Lambda_{l\;k}(g)=(-1)^{k+l}\Lambda_{n-l\;n-k}(g)$.
The action of
$\mathfrak{su}(2)$ on $V_n$ is as follows:
\begin{equation}\label{2.3}
\textbf{v}_{k,n}d\rho_{n}(\varepsilon):=\frac{d}{dt}
\big(\rho_n(\exp{t\varepsilon})(\textbf{v}_{k,n})\big)\big|_{t=0},
\end{equation}
for $0\leq k\leq n$ and any element
$\varepsilon\in\mathfrak{su}(2)$. In particular, when $\varepsilon$ takes
$\varepsilon_1$, $\varepsilon_2$, $\varepsilon_3$ respectively, we
have
\begin{equation}\label{2.4}
\textbf{v}_{k,n}d\rho_{n}(\varepsilon_1)=(n-2k)\texttt{i}\textbf{v}_{k,n},
\end{equation}
\begin{equation}\label{2.5}
\textbf{v}_{k,n}d\rho_{n}(\varepsilon_2)=-a_{k-1,n}\textbf{v}_{k-1,n}
+a_{k,n}\textbf{v}_{k+1,n},
\end{equation}
\begin{equation}\label{2.6}
\textbf{v}_{k,n}d\rho_{n}(\varepsilon_3)=
-a_{k-1,n}\texttt{i}\textbf{v}_{k-1,n}-a_{k,n}\texttt{i}
\textbf{v}_{k+1,n},
\end{equation}
for $0\leq k\leq n$ and $a_{k,n}=\sqrt{(k+1)(n-k)}$. Set
$$\lambda_k:=n-2k,\;\;\;\;\;\textbf{v}_{\lambda_k,n}:=\textbf{v}_{k,n},\;\;\;\;
\;a_{\lambda_k,n}=a_{k,n}=\sqrt{[(n+1)^2-(\lambda_k-1)^2]}\big/2.$$
An element in $\Delta_n=\{\lambda_0,\ldots,\lambda_n\}$ is called a weight of $\rho_n$, and $\lambda_0$ is called the highest weight. The representation space $V_n$ has an orthonormal decomposition w.r.t. the Hermitian inner product $(\;,\;)$, i.e., $V_n=V_{\lambda_0}\oplus\cdots\oplus V_{\lambda_n}$, where $V_{\lambda_k}=\mbox{Span}_{\textbf{C}}\{\textbf{v}_{\lambda_k,n}\}$ is called the weight space w.r.t. the weight $\lambda_k$.

Suppose $n=2m-1\in\textbf{N}^+$ is an odd number. Define $\textbf{u}_{\lambda_k,n}$, $\textbf{u}_{-\lambda_k,n}$ as follows
\begin{equation}\label{2.7}
\textbf{u}_{\lambda_k,n}:=\textbf{v}_{\lambda_k,n},\;\;\;
\textbf{u}_{-\lambda_k,n}:=(-1)^k\texttt{i}\textbf{v}_{-\lambda_k,n},
\end{equation}
where $0\leq k\leq m-1$. Under the basis $\{\textbf{u}_{\lambda_0,n},\ldots,\textbf{u}_{\lambda_{m-1},n},
\textbf{u}_{-\lambda_0,n},\ldots,
\textbf{u}_{-\lambda_{m-1},n}\}$, for every $g\in SU(2)$, one can check that $\rho_n(g)$ satisfies
\begin{equation}\label{2.8}
J\;\rho_n(g)=\overline{\rho_n(g)}\;J,\;\;\;\mbox{for}\;\;\;J=\left(
                                               \begin{array}{cc}
                                                 0 & -I_m \\
                                                 I_m & 0 \\
                                               \end{array}
                                             \right).
\end{equation}
We identify $V_{2m-1}$ with $\textbf{H}^m$ by
\begin{equation}\label{2.9}
\textbf{v}=\sum\limits_{k=0}^{m-1}(a_k\textbf{u}_{\lambda_k,n}+b_k
\textbf{u}_{-\lambda_k,n})\mapsto\big(a_0+b_0\texttt{j},\ldots,a_{m-1}+
b_{m-1}\texttt{j}\big),
\end{equation}
where $\texttt{j}\in\textbf{H}$ and $\texttt{j}^2=-1$. From this identification, it is convenient to define
\begin{equation}\label{2.10}
{\texttt{j}\textbf{u}_{\lambda_k,n}:=\textbf{u}_{-\lambda_k,n}}.
\end{equation}
 Thus, an element in $V_{2m-1}$ can be written as $\textbf{v}=\sum\limits_{k=0}^{{m-1}}(a_k+b_k{
 \texttt{j}}){\textbf{u}_{\lambda_k,n}}$. When $V_{2m-1}$ is viewed as a quaternionic linear space, the property (\ref{2.8}) ensures us obtain a quaternionic representation of $SU(2)$ denoted by $\rho_m$. In terms of matrix, we have
\begin{equation}\label{2.11}
\rho_m:SU(2)\longrightarrow Sp(m),\;\;\;g\mapsto \Big(\Xi_{k\;l}(g):=\Lambda_{k\;l}(g)+(-1)^{k+1}\texttt{i}\Lambda_{k\;2m-1-l}(g)\texttt{j}\Big),
\end{equation}
where $0\leq k,l\leq m-1$ and $Sp(m)=\big\{A\in M(m,\textbf{H}):AA^*=I_m\big\}$ is the unitary symplectic group.

It's well known that $\big\{(V_n,\rho_n):n=0,1,2,\cdots\big\}$ are all inequivalent unitary representations of $SU(2)$. By the Theorem 6.3 in \cite{[BD]}, we know that $\big\{(V_{2m-1},\rho_m):m=1,2,\ldots\big\}$ are all inequivalent proper unitary quaternionic representations of $SU(2)$.
Since any quaternionic representation $(V,\rho)$ is completely reducible, for a proper one, we can write
\begin{equation}\label{2.12}
V=\bigoplus\limits_{\alpha=1}^s\;V_{2m_\alpha-1},\;\;\;\rho=\bigoplus
\limits_{\alpha=1}^s\;\rho_{m_\alpha}.
\end{equation}
In terms of weights, under the natural inner product, we also have the orthogonal decomposition
\begin{equation}\label{2.13}
V=\bigoplus\limits_{\lambda>0}\;U_\lambda=\bigoplus\limits_{\lambda>0}
\Big(\bigoplus\limits_\alpha\;U_{\lambda,\alpha}\Big),\;\;\;U_{\lambda,\alpha}
=V_{\lambda,2m_\alpha-1}\oplus V_{-\lambda,2m_\alpha-1}.
\end{equation}
Here, the "orthogonal decomposition" is related to the inner product $(\;,\;)$ defined by
\begin{equation*}
(\textbf{z},\textbf{w}):=\sum\limits_{k=0}^n\;p_k\bar q_k,
\end{equation*}
for $\textbf{z}=\sum\limits_{k=0}^np_k\textbf{u}_{\lambda_k,n}$,
$\textbf{w}=\sum\limits_{k=0}^nq_k\textbf{u}_{\lambda_k,n}$.
Notice that $d\rho$ is a real representation of $\mathfrak{su}(2)$, we extend naturally $d\rho$ to be a complex representation of $\mathfrak{sl}(2;\textbf{C})$, which is still denoted by $d\rho$. Set
$$\sigma_1=\left(
             \begin{array}{cc}
               1 & 0 \\
               0 & -1 \\
             \end{array}
           \right),\;\;\;\sigma_2=\left(
                                    \begin{array}{cc}
                                      0 & 1 \\
                                      0 & 0 \\
                                    \end{array}
                                  \right),\;\;\;
                                  \sigma_3=\left(
                      \begin{array}{cc}
                        0 & 0 \\
                        1 & 0 \\
                      \end{array}
                    \right),$$
which forms a basis of $\mathfrak{sl}(2;\textbf{C})$. In the sequel, we denote $d\rho(\sigma_1),d\rho(\sigma_2),d\rho(\sigma_3)\in End(V)$ by $H,A,B$ respectively. Alternatively, they are determined by
\begin{equation}\label{2.14}
H=-\texttt{i}d\rho(\varepsilon_1),\;\;\;
A=\frac{d\rho(\varepsilon_2)-\texttt{i}d\rho(\varepsilon_3)}{2},\;\;\;
B=-\frac{d\rho(\varepsilon_2)+\texttt{i}d\rho(\varepsilon_3)}{2},
\end{equation}
which satisfy
\begin{equation}\label{2.15}
[H,A]=2A,\;\;\;[H,B]=-2B,\;\;\;[A,B]=H.
\end{equation}
It follows from (\ref{2.4})-(\ref{2.7}), (\ref{2.10}) and (\ref{2.14}), we have
\begin{eqnarray}\label{2.16}
\textbf{u}_{\lambda_k,n_\alpha}H=\lambda_k\textbf{u}_{\lambda_k,n_\alpha},\;\;\;
(\texttt{j}\textbf{u}_{\lambda_k,n_\alpha})H=-\lambda_k(\texttt{j}
\textbf{u}_{\lambda_k,n_\alpha}),\;\;\mbox{for}\;\;0\leq k\leq m_\alpha-1,\nonumber\\
\textbf{u}_{\lambda_k,n_\alpha}A=-a_{\lambda_k,n_\alpha}\textbf{u}_{\lambda_k-2,
n_\alpha},\;\;\mbox{for}\;\;0\leq k\leq m_\alpha-2,\nonumber\\
{\textbf{u}_{\lambda_{m_\alpha-1},n_\alpha}A=-(-1)^{m_\alpha}\texttt{i}
a_{\lambda_{m_\alpha-1},n_\alpha}(\texttt{j}\textbf{u}_{\lambda_{m_\alpha-1},n_\alpha})},
\hspace{3.3cm}\nonumber\\
(\texttt{j}\textbf{u}_{\lambda_k,n_\alpha})A=a_{-\lambda_k,n_\alpha}(\texttt{j}
\textbf{u}_{\lambda_k,n_\alpha}),\;\;\mbox{for}\;\;0\leq k\leq m_\alpha-1,\\
\textbf{u}_{\lambda_k,n_\alpha}B=-a_{\lambda_k+2,n_\alpha}\textbf{u}_{\lambda_k+2,
n_\alpha},\;\;\mbox{for}\;\;0\leq k\leq m_\alpha-1,\nonumber\\
(\texttt{j}\textbf{u}_{\lambda_k,n_\alpha})B=a_{-\lambda_k+2,n_\alpha}
(\texttt{j}\textbf{u}_{\lambda_k-2,n_\alpha}),\;\;\mbox{for}\;\;0\leq k \leq m_\alpha-2,\nonumber\\
{(\texttt{j}\textbf{u}_{\lambda_{m_\alpha-1},n_\alpha})B=(-1)^{m_\alpha}
\texttt{i}a_{\lambda_{m_\alpha-1},n_\alpha} \textbf{u}_{\lambda_{m_\alpha-1},n_\alpha}}\nonumber.\hspace{3.3cm}
\end{eqnarray}
The pull-back of Maurer-Cartan forms of $Aut(V)$ is
\begin{equation}\label{2.17}
d\rho\rho^{-1}=(\texttt{i}\omega)H+\varphi A-\bar{\varphi}B.
\end{equation}
It follows that
\begin{equation}\label{2.-17}
\textbf{z}d\rho=\big((\texttt{i}\omega)(\textbf{z}H)+\varphi(\textbf{z} A)-\bar{\varphi}(\textbf{z}B)\big)\rho.
\end{equation}

Throughout this paper we will agree on the following conventions:
\begin{itemize}
\item $n_\alpha=2m_\alpha-1$, for $m_\alpha\in\textbf{N}^+$;

    \item
    $\Delta_\alpha:=\{n_\alpha,n_\alpha-2,\cdots,-(n_\alpha\!-\!2),-n_\alpha\}$;

    \item$a_{\lambda,\alpha}:=a_{\lambda,n_\alpha}$, and
$a_{\lambda,\alpha}=0$ if $\lambda\notin\Delta_\alpha$;

\item $\textbf{u}_{\lambda,\alpha}:=\textbf{u}_{\lambda,n_\alpha}$,
$V_{\lambda,\alpha}=\mbox{Span}_{\mathbb{C}}\{\textbf{u}_{\lambda,\alpha}\}$,
and $V_{\lambda,\alpha}=\{0\}$ if $\lambda\notin\Delta_\alpha$;

\item {$\phi_{\lambda,m}:=\textbf{u}_{\lambda,n}\rho_m(g)$, $g\in SU(2)$, for $\lambda>0$}.

  \end{itemize}

\subsection{Geometry of surfaces in quaternionic projective space $\textbf{HP}^n$}
Throughout this subsection we agree on the following ranges of indices:
\begin{equation*}
0\leq A,B,\ldots\leq n,\;\; 1\leq \alpha,\beta,\ldots\leq n.
\end{equation*}
The standard inner product $(\;,\;)$ of $\textbf{H}^{n+1}$ is defined by
\begin{equation*}
(\textbf{z},\textbf{w})=\sum\limits_{A}z_A\bar w_A=\sum\limits_{A}
z_Aw_{\bar A},
\end{equation*}
for $\textbf{z}=(z_0,z_1,\ldots,z_n)$, $\textbf{w}=(w_0,w_1,\ldots,w_n)\in\textbf{H}^{n+1}$. A unitary frame of $\textbf{H}^{n+1}$ is an ordered set of $n+1$ linearly independent vectors $Z_0,Z_1,\ldots,Z_n$ satisfying
\begin{equation*}
(Z_A,Z_B)=\delta_{AB}.
\end{equation*}
Taking the exterior derivative of $Z_A$, we have
\begin{equation}\label{2.18}
dZ_A=\sum\limits_{B}\omega_{A\bar B}\,Z_B,
\end{equation}
where $\omega_{A\bar B}$ are $\textbf{H}$-valued one forms. We identify the space of unitary frames with the unitary symplectic group $Sp(n+1)$. Then, $\omega_{A\bar B}$ are the Maurer-Cartan forms of $Sp(n+1)$, which satisfy the Maurer-Cartan equations:
\begin{equation}\label{2.19}
d\omega_{A\bar B}=\sum\limits_{C}\omega_{A\bar C}\wedge\omega_{C\bar B}.
\end{equation}

The quaterionic projective space $\textbf{HP}^n$ is the set of all one-dimensional subspaces in $\textbf{H}^{n+1}$. An element of $\textbf{HP}^n$ can be defined by the unitary vector $Z_0$,  up to a factor of absolute 1. Its orthogonal vectors $Z_\alpha$ are defined up to a transformation of $Sp(n)$. So, $\textbf{HP}^{n}$ has a $Sp(1)\times Sp(n)$-structure. Therefore, the form
\begin{equation}\label{2.20}
ds^2=\sum\limits_{\alpha}\omega_{0\bar\alpha}\,\omega_{\bar 0\alpha}
\end{equation}
defines a Riemannian metric on $\textbf{HP}^{n}$.

Let $M$ be an oriented Riemnnian surface and let $f:M\longrightarrow
\textbf{HP}^{n}$ be a smooth immersion. The induced metric on $M$ by
$ds^2_M=\varphi\bar\varphi$, where $\varphi$ is complex-valued one-form. $\varphi$ is defined up to a complex factor of absolute value 1.
The first structure equation of $ds^2_M$ is given by
\begin{equation}\label{2.21}
d\varphi=\eta\wedge\varphi,
\end{equation}
where the one-form $\eta$ is purely imaginary. The Gauss curvature $K$ is described by
\begin{equation}\label{2.22}
d\eta=-\frac{K}{2}\varphi\wedge\bar{\varphi}.
\end{equation}
Locally, we have
\begin{equation}\label{2.23}
f^*\omega_{0\bar{\alpha}}=\varphi\,a_{\alpha}+\bar\varphi\,
b_{\alpha},
\end{equation}
where $a_{\alpha}$, $b_{\alpha}$ are local $\textbf{H}$-valued functions. Taking the exterior derivative on both sides of (\ref{2.23}) and using (\ref{2.19}), (\ref{2.21}) and (\ref{2.22}), we obtain
\begin{equation}\label{2.24}
\varphi\wedge Da_{\alpha}+\bar{\varphi}\wedge Db_{\alpha}=0,
\end{equation}
where
\begin{equation}\label{2.25}
Da_{\alpha}=da_{\alpha}+\eta a_{\alpha}-\theta_{0\bar 0}a_{\alpha}-\theta'_{0\bar 0}b_{\alpha}+a_{\beta}\omega_{\beta\bar\alpha},
\end{equation}
\begin{equation}\label{2.26}
Db_{\alpha}=db_{\alpha}-\eta b_{\alpha}-\theta'_{0\bar 0}a_{\alpha}-\theta_{0\bar 0}b_{\alpha}+b_{\beta}\omega_{\beta\bar\alpha},
\end{equation}
and
\begin{equation}\label{2.27}
\theta_{0\bar 0}=\frac{\omega_{0\bar 0}-\texttt{i}\omega_{0\bar 0}\texttt{i}}{2},\;\;\;\theta'_{0\bar 0}=\frac{\omega_{0\bar 0}+\texttt{i}\omega_{0\bar 0}\texttt{i}}{2}.
\end{equation}
Because $\textbf{HP}^{n}$ is non-Kahlerian, the definition
of $Da_{\alpha}$ (resp. $Db_{\alpha}$) involves $b_{\alpha}$ (resp. $a_{\alpha}$). By (\ref{2.24}), we can set
\begin{equation}\label{2.28}
Da_{\alpha}=\varphi\,p_{\alpha}+\bar{\varphi}\,q_{\alpha},
\;\;\;\;\;
Db_{\alpha}=\varphi\,q_{\alpha}+\bar{\varphi}\,r_{\alpha}.
\end{equation}
Similar to \cite{[CW2]}, $f$ is minimal if and only if
$q_{\alpha}=0$, i.e.,
\begin{equation}\label{2.29}
Da_{\alpha}\equiv0,\mod\varphi\;\;\;\;\text{or}\;\;\;\;
Db_{\alpha}\equiv0,\mod\bar\varphi.
\end{equation}
Set the local $\textbf{H}^{n+1}$-valued functions $X,Y$ as
\begin{equation}\label{2.30}
X=\sum\limits_{\alpha}a_{\alpha}Z_\alpha,\;\;\;\;
Y=\sum\limits_{\alpha}b_{\alpha}Z_\alpha,
\end{equation}
and define the covariant derivative of $X$ and $Y$ by
\begin{equation}\label{2.31}
DX=dX+\eta X-\theta_{0\bar 0}X-\theta'_{0\bar 0}Y,
\end{equation}
\begin{equation}\label{2.32}
DY=dY-\eta Y-\theta_{0\bar 0}X-\theta'_{0\bar 0}Y.
\end{equation}
We give a criterion to measure the minimality of $f$ in terms of $DX$ or $DY$, that is
\begin{prop}\label{proposition2.1}
The smooth immersion $f$ is minimal if and only if one of the following holds:

\emph{(a)} $DX\equiv0,\;\;\;\mod Z_0,\;\varphi;$

\emph{(b)} $DY\equiv0,\;\;\;\mod Z_0,\;\bar\varphi.$
\end{prop}

\emph{Proof}.
It follows from the facts
\begin{equation*}
DX=\sum\limits_{\alpha}Da_{\alpha}\;Z_\alpha+\sum\limits_{\beta}
a_{\beta}\omega_{\beta\bar 0}\;Z_0,\;\;\;
DY=\sum\limits_{\alpha}Db_{\alpha}\;Z_\alpha+\sum\limits_{\beta}
b_{\beta}\omega_{\beta\bar 0}\;Z_0.
\end{equation*}
This completes the proof. \hfill$\Box$

\section{Homogeneous two-spheres in $\textbf{HP}^n$}\label{sec3}
In this section we first determine the relations between the homogeneous two-spheres and the two-dimensional $\rho(SU(2))$-orbits in $\textbf{HP}^n$, where $\rho$ is a quaternionic representation of $SU(2)$. Then, all the two-dimensional
$\rho(SU(2))$-orbits in $\textbf{HP}^n$ are determined.

\begin{theo}\label{theorem 3.1}
Let $f:S^2\longrightarrow \textbf{\emph{HP}}^n$ be a homogeneous immersion, then
there exists a quaternionic representation $\rho$ of $SU(2)$ such that $f(S^2)$ is a
two-dimensional $\rho(SU(2))$-orbit in $\textbf{\emph{HP}}^n$.
\end{theo}

\emph{Proof}. It is similar to the proof of Theorem 3.2 in \cite{[FJXX]} or Theorem 3.3 in \cite{[PJX]}, we
omit the details.\hfill$\Box$

\vspace{0.3cm}

In general, a $\rho(SU(2))$-orbit in $\textbf{HP}^n$ is a principle orbit,
i.e., it is a three-dimensional orbit. For the two-dimensional
$\rho(SU(2))$-orbits, we have
\begin{lem}\label{lemma3.2}
Let $M$ be a $\rho(SU(2))$-orbit in $ \textbf{\emph{HP}}^n$, then $M$ is
two-dimensional if and only if there exists a point $P_0\in M$ s.t.
$P_0$ is invariant under the $\rho(T)$-action.
\end{lem}

\emph{Proof}. Let $P$ be a point in $M$. Since $dim M=2$, the isotropy group
$G$ at $P$ is a 1-dimensional subgroup of $SU(2)$. Therefore, there exists a nonzero
element $\varepsilon\in\mathfrak{su}(2)$ s.t.
$T_{I_2}G=\mbox{Span}_{\mathbb{R}}\{\varepsilon\}$. Notice that
$\varepsilon\in\mathfrak{su}(2)$, there exist a fixed $g\in
SU(2)$ s.t. $g^{-1}\varepsilon g= c \varepsilon_1$ for some
nonzero real number $c$. It follows that the isotropy group at
$P_0=P\rho(g)\in M$ contains the maximal torus subgroup $T$,
i.e., $P_0$ is invariant under $\rho(T)$-action. The sufficiency
follows from the fact that there is no two-dimensional subgroup in $SU(2)$.
\hfill$\Box$

\vspace{0.3cm}

\noindent\textbf{\emph{Remark}}. By taking derivative, the condition "$P_0$ is invariant under $\rho(T)$" is equivalent to
\begin{equation}\label{3.1}
H(P_0)\subseteq P_0,
\end{equation}
where $P_0$ is viewed as a linear subspace and $H$ is defined in (\ref{2.14}).

\vspace{0.3cm}

Let $\rho$ be a quaternonic representation of $SU(2)$. Then $\rho(SU(2))$ acts on $\textbf{HP}^n$ by
\begin{equation}\label{3.2}
SU(2)\times \textbf{HP}^n\longrightarrow \textbf{HP}^n,\;(g,P)\mapsto P\rho(g).
\end{equation}
Therefore, for a two-dimensional $\rho(SU(2))$-orbit
$M$ in $\textbf{HP}^n$, we can obtain an immersion $f$ from $S^2$ into
$\textbf{HP}^n$ as follows
\begin{equation}\label{3.3}
f:S^2\simeq SU(2)/T\longrightarrow \textbf{HP}^n,\; [g]\mapsto P_0\rho(g),
\end{equation}
where $P_0$ is a fixed point in $M$. According to Theorem \ref{theorem
3.1}, we only need to study the two-dimensional $\rho(SU(2))$-orbits
in $\textbf{HP}^n$ if we study the homogeneous immersions from
$S^2$ into $\textbf{HP}^n$.

\vspace{0.3cm}

Let $\rho=\rho_{m_1}\oplus\cdots\oplus\rho_{m_s}$ be a quaternionic
representation of $SU(2)$. Let $M$ be a linearly full two-dimensional $\rho(SU(2))$-orbit in $\textbf{HP}^n$. Suppose $P_0=[\textbf{z}]$ is a point obtained in
Lemma \ref{lemma3.2}. According to the decomposition
(\ref{2.13}), we can write
$$\textbf{z}=\sum\limits_{\lambda>0}\sum\limits_{\alpha}(c_{\lambda,\alpha}
+d_{\lambda,\alpha}\texttt{j})\textbf{u}_{\lambda,\alpha},\;\;\mbox{for}\;\;
c_{\lambda,\alpha},d_{\lambda,\alpha}\in\textbf{C}.$$
Without loss of generality, we further assume $d_{\lambda_0,\alpha_0}=0$ for some $\lambda_0$ and $\alpha_0$. If not, we use $(\bar{c}_{\lambda_0,\alpha_0}-d_{\lambda_0,\alpha_0}\texttt{j})\textbf{z}$ instead of $\textbf{z}$. The condition (\ref{3.1}) reduces to
\begin{equation}\label{3.4}
\textbf{z}H=p\textbf{z},
\end{equation}
for some fixed $p=p_1+p_2\texttt{j}\in \textbf{H}$. In terms of coefficients, (\ref{3.4}) is equivalent to
\begin{equation}\label{3.5}
\lambda c_{\lambda,\alpha}=p_1c_{\lambda,\alpha}-p_2\bar{d}_{\lambda,\alpha},\;\;\;
-\lambda d_{\lambda,\alpha}=p_1d_{\lambda,\alpha}+p_2\bar{c}_{\lambda,\alpha},
\end{equation}
for all $\lambda,\alpha$. Notice that $d_{\lambda_0,\alpha_0}=0$, it follows from (\ref{3.5}) we have $p_1=\lambda_0$, $p_2=0$. This implies
\begin{equation}\label{3.6}
\lambda=p_1=\lambda_0,\;\;
d_{\lambda,\alpha}=0, \;\;\mbox{for all}\;\; \lambda,\alpha,
\end{equation}
by the linearly full assumption. Thus, we have
\begin{equation}\label{3.7}
\textbf{z}=\sum\limits_{\alpha}c_\alpha\textbf{u}_{\lambda,\alpha},
\;\;\;c_\alpha\in\textbf{C},
\end{equation}
for a positive weight $\lambda$.

In summery, we completely classify all homogeneous two-spheres in $\textbf{HP}^{n}$. That is
\begin{theo}\label{theorem3.3}
Let $\rho=\rho_{m_1}\oplus\cdots\oplus\rho_{m_s}$, $m_1+\cdots+m_s=n+1$, be a quaternionic representation of $SU(2)$. If $M$ is a two-dimensional $\rho(SU(2))$-orbit in $\textbf{\emph{HP}}^{n}$, then there exists a point $[\textbf{\emph{z}}]\in M$ and a positive weight $\lambda$ s.t.
$$\textbf{\emph{z}}=\sum\limits_{\alpha}c_\alpha
\textbf{\emph{u}}_{\lambda,\alpha},\;\;\;c_\alpha\in\textbf{\emph{C}}.$$
%\emph{(2)} If $M$ is a two-dimensional $\rho(SU(2))$-orbit in $G(2,n;\textbf{\emph{H}})$, there exists a point $P=\mbox{\emph{Span}}_{\textbf{\emph{H}}}\{\textbf{\emph{z}},\textbf{\emph{w}}\}$
%and positive weights $\lambda,\mu$ s.t.
%$$\textbf{z}=\sum\limits_{\alpha}c_\alpha\textbf{\emph{u}}_{\lambda,\alpha},\;\;\;
%\textbf{w}=\sum\limits_{\alpha}c'_\alpha\textbf{\emph{u}}_{\mu,\alpha},
%\;\;\;c_\alpha,\;c'_{\alpha}\in\textbf{\emph{C}}.$$
\end{theo}

\vspace{0.3cm}

%\noindent\textbf{\emph{Remark}.} This theorem can be generalized to a more general setting. Namely, for each
%two-dimensional $\rho(SU(2))$-orbit $M$ in $G(k,n;\textbf{H})$, there exists a
%point $P_0=\mbox{Span}_{\textbf{H}}\{\textbf{z}_1,\ldots,\textbf{z}_k\}\in M$ s.t. $\textbf{z}_i=\sum\limits_{\alpha}c_{i,\alpha}\textbf{u}_{\lambda_i,\alpha}$ for some weights $\lambda_1,\ldots,\lambda_k$.

\section{Minimal homogeneous two-spheres in $\textbf{HP}^n$}\label{sec4}

Let $\rho=\rho_{m_1}\oplus\cdots\oplus\rho_{m_s}$, $m_1+\cdots+m_s=n+1$ be a quaternionic representation of $SU(2)$, and let $M$ be a two-dimensional $\rho(SU(2))$-orbit in $\textbf{HP}^n$. Suppose $[\textbf{z}]\in M$ is the point we get in Theorem \ref{theorem3.3}, i.e., $\textbf{z}=\sum\limits_{\alpha}\;c_\alpha\textbf{u}_{\lambda,\alpha}$, $c_\alpha\in\textbf{C}$, for some positive weight $\lambda$. We further assume $|\textbf{z}|^2=1$. Set
\begin{equation}\label{4.1}
Z_0=\textbf{z}\rho.
\end{equation}
Taking the exterior derivative of (\ref{4.1}) and using (\ref{2.16}), (\ref{2.17}), we have
\begin{eqnarray}\label{4.2}
dZ_0=\big(\lambda(\texttt{i}\omega)\textbf{z}+\varphi\textbf{z}A-\bar{\varphi}\textbf{z}B\big)\rho.
\end{eqnarray}
We define
\begin{equation}\label{4.3}
dZ_0\equiv \varphi X+\bar{\varphi} Y,\;\;\; \mod Z_0.
\end{equation}
Notice that $\textbf{z}\in U_{\lambda}$, $\textbf{z}A\in U_{\lambda-2}$ for $\lambda>1$, $\textbf{z}A\in U_{\lambda}$ for $\lambda=1$,  $\textbf{z}B\in U_{\lambda+2}$, we have
\begin{equation}\label{4.4}
X=\textbf{z}A\rho-\ell\textbf{z}\rho,\;\;\;
Y=-\textbf{z}B\rho,
\end{equation}
where $\ell=(\textbf{z}A,\textbf{z})$. It's clear that $\ell=0$ when $\lambda>1$, and $\ell$ has the possibility to be nonzero when $\lambda=1$.
The fact $(X,Y)=0$ means that the orbit $M$ is conformal in $\textbf{HP}^n$, and the induced metric is
\begin{equation}\label{4.5}
(|X|^2+|Y|^2)\;\varphi\bar\varphi.
\end{equation}
Notice that $|X|^2+|Y|^2$ is a constant and the metric $\varphi\bar{\varphi}$ has the curvature 4, we know that the induced metric (\ref{4.5}) has constant curvature
\begin{equation}\label{4.6}
K=\frac{4}{|X|^2+|Y|^2}.
\end{equation}

Now we give a simple expression of the criteria (a) in Proposition \ref{proposition2.1} for later use. Notice that
$\eta=2\texttt{i}\omega$ by (\ref{2.2}), and substituting (\ref{4.2}) into (\ref{2.27}), we have $\theta_{0\bar 0}=\texttt{i}\omega$, $\theta'_{0\bar 0}=\varphi\ell$. From the definition (\ref{2.31}), we get
\begin{equation}\label{4.7}
DX=\varphi\Big(\textbf{z}A^2-(\ell\textbf{z})A-(\ell\textbf{z})B)\Big)\rho
-\bar\varphi\Big(\textbf{z}AB-(\ell\textbf{z})B\Big)\rho.
\end{equation}
Thus, according to Proposition \ref{proposition2.1}, $M$ is minimal in
$\textbf{HP}^{n}$ if and only if
\begin{equation}\label{4.8}
\textbf{z}AB-(\ell\textbf{z})B=p\textbf{z}
\end{equation}
for fixed $p=p_1+p_2\texttt{j}\in\textbf{H}$. We will study the geometry of  minimal orbits in following two cases: $\lambda>1$ and $\lambda=1$.

\vspace{0.3cm}

\noindent\textbf{Case I.} \hspace{0.3cm} Suppose that $\lambda>1$, and hence $\ell=0$. From (\ref{2.16}), by comparing the coefficients, (\ref{4.7}) is equivalent to
\begin{equation}\label{4.9}
a_{\lambda,\alpha}^2=p,\;\;\;\mbox{for all}\;\; \alpha.
\end{equation}
This implies that $n_1=\cdots=n_s=:n_0=2m_0-1$. So, the associated quaternionic representation $\rho$ takes the form $\rho=\rho_{m_0}\oplus\cdots\oplus\rho_{m_0}$, and the corresponding minimal immersion is given by
\begin{equation}\label{4.10}
[a,b]\mapsto[c_1\phi_{\lambda,m_0},\ldots,c_s\phi_{\lambda,m_0}],
\end{equation}
where $c_\alpha\in\textbf{C}$ are nonzero.
\begin{prop}\label{proposition4.1}
Up to an isometry of $\textbf{\emph{HP}}^{n}$, the immersion \emph{(\ref{4.10})} is congruent to
\begin{equation}\label{4.11}
f_{\lambda,m_0,\;t}:[a,b]\mapsto [\cos t\;\phi_{\lambda,m_0},\texttt{\emph{i}}\sin t\;\phi_{\lambda,m_0},0,\ldots,0],
\end{equation}
where $t\in [0,\pi/2).$
\end{prop}

\emph{\textbf{Remark}}. If the parameter $t=0$ and $f$ is proper, we obtain the immersion
\begin{equation*}
f_{\lambda}:S^2\longrightarrow \textbf{HP}^n,\hspace{0.5cm} [a,b]\mapsto[\phi_{\lambda,n+1}],
\end{equation*}
which firstly appeared in \cite{[Oh]}.

\vspace{0.3cm}

\emph{Proof}. Set $\textbf{c}:=(c_1,\ldots,c_s)$. Notice that the choice of $\textbf{c}$ is unique up to a factor (multiplied on the left) of absolute 1, without loss of generality, we can find a matrix $\textbf{T}\in SO(s)$ s.t. $\textbf{c}\;\textbf{T}=(\cos t,\texttt{i}\sin t,0,\ldots,0)$. Then the proposition holds by choosing the isometry $\textbf{T}\otimes I_{m_0}$ of $\textbf{HP}^{n}$. \hfill$\Box$

\vspace{0.3cm}
Next, we calculate the Gauss curvature of the immersion (\ref{4.11}). By (\ref{2.16}) and (\ref{4.4}), we obtain
$$|X|^2+|Y|^2=|\textbf{z}A|^2+|\textbf{z}B|^2=a^2_{\lambda,n_0}+a^2_{\lambda+2,n_0}.$$
So by (\ref{4.6}), the Gauss curvature is
\begin{equation}\label{4.12}
K=\frac{4}{a^2_{\lambda,n_0}+a^2_{\lambda+2,n_0}}=\frac{8}{(n_0+1)^2-(\lambda^2+1)}.
\end{equation}
\vspace{0.3cm}

\noindent\textbf{Case II.} \hspace{0.3cm}Suppose that $\lambda=1$. Since $\textbf{z}=\sum\limits_{\alpha=1}^sc_\alpha
\textbf{u}_{1,\alpha}$, $c_\alpha\in\textbf{C}$, we have $$\textbf{z}A=\sum\limits_{\alpha=1}^s\big((-1)^{m_\alpha+1}\texttt{i}c_{\alpha}a_{1,\alpha}
\texttt{j}\big)\textbf{u}_{1,\alpha}$$ by (\ref{2.16}), and then,
\begin{equation}\label{4.13}
(\textbf{z}A,\textbf{z})=\Big(\sum\limits_{\alpha=1}^s(-1)^{m_\alpha+1}
a_{1,\alpha}c_\alpha^2\texttt{i}\Big)\texttt{j}.
\end{equation}
From (\ref{4.13}), we find $\big((e^{\texttt{i}\tau}\textbf{z})A,e^{\texttt{i}\tau}\textbf{z}\big)=
e^{2\texttt{i}\tau}(\textbf{z}A,\textbf{z})$. So, we can further assume
$\ell=(\textbf{z}A,\textbf{z})=(\ell'\texttt{i})\texttt{j}$, where $\ell'\in\textbf{R}$. In terms of components, (\ref{4.8}) is equivalent to
\begin{equation}\label{4.14}
a_{1,\alpha}^2c_\alpha+(-1)^{m_\alpha}\ell'a_{1,\alpha}\bar c_\alpha
=p_1c_\alpha,\;\;\;\mbox{for all}\;\;\alpha.
\end{equation}
Splitting $c_\alpha$, $p_1$ into real and imaginary parts, i.e., $c_\alpha=
c'_\alpha+c''_{\alpha}\texttt{i}$, $p_1=p'_1+p''_1\texttt{i}$. Then, (\ref{4.14}) are equivalent to
\begin{equation}\label{4.15}
(a_{1,\alpha}^2+(-1)^{m_\alpha}\ell'a_{1,\alpha}-p_1')c_\alpha'+p_1''c_\alpha''=0,
\end{equation}
\begin{equation}\label{4.16}
-p_1''c_\alpha'+(a_{1,\alpha}^2-(-1)^{m_\alpha}\ell'a_{1,\alpha}-p_1')c_\alpha''=0,
\end{equation}
for all $\alpha$. Notice that $c_\alpha\neq 0$, we know that the determinant of the coefficient matrix of the system (\ref{4.15}) and
(\ref{4.16}) is equal to zero, which follows that $a_{1,\alpha}^2$ satisfy the real quadratic equation
\begin{equation*}
x^2-\big((\ell')^2+2p_1'\big)x+|p_1|^2=0.
\end{equation*}
This means that there are two different $m_\alpha$ at most. So, the associated representation $\rho$ can be written as
\begin{equation}\label{4.17}
\rho=\rho_{m_1}\oplus\cdots\oplus\rho_{m_1},
\end{equation}
or
\begin{equation}\label{4.18}
\rho=\rho_{m_1}\oplus\cdots\oplus\rho_{m_1}\oplus\rho_{m_2}
\oplus\cdots\oplus\rho_{m_2},\;\;\;m_1<m_2.
\end{equation}

\textbf{(II.I)} \hspace{0.3cm}If the associate representation $\rho$ takes form of (\ref{4.17}), in terms of homogeneous coordinates, the corresponding immersion is given by
\begin{equation}\label{4.19}
[a,b]\mapsto[c_1\phi_{1,m_1},\ldots,c_s\phi_{1,m_1}],
\end{equation}
where $c_\alpha\in\textbf{C}$.

\begin{prop}\label{proposition4.2}
Up to an isometry of $\textbf{\emph{HP}}^{n}$, the minimal immersion \emph{(\ref{4.19})}
is congruent to
\begin{equation}\label{4.20}
f_1:[a,b]\mapsto[\phi_{1,m_1},0,\ldots,0]
\end{equation}
or
\begin{equation}\label{4.21}
f_{m_1,\;m_1}:[a,b]\mapsto\frac{\sqrt{2}}{2}[\phi_{1,m_1},\texttt{\emph{i}}\;\phi_{1,m_1},
0,\ldots,0].
\end{equation}
\end{prop}

\emph{Proof}. For this subcase, we have $\ell'=(-1)^{m_1+1}a_{1,n_1}\sum\limits_{\alpha}c_{\alpha}^2$. Then the minimality condition (\ref{4.14}) can be written as
\begin{equation}\label{4.22}
(a_{1,n_1}^2-p_1)c_\alpha=a_{1,n_1}^2(\sum\limits_{\alpha}c_{\alpha}^2)\bar{c}_\alpha
\end{equation}
for all $\alpha$, where $\sum\limits_{\alpha}c_{\alpha}\bar{c}_\alpha=1$. Multiplying by $c_{\alpha}$ on both sides of (\ref{4.22}) and summating for all $\alpha$, we obtain $$(a_{1,n_1}^2-p_1)\sum\limits_{\alpha}c_{\alpha}^2=a_{1,n_1}^2(\sum\limits_{\alpha}c_{\alpha}^2),$$
which implies $p_1=0$ or $\sum\limits_{\alpha}c_{\alpha}^2=0$.

If $p_1=0$, then from (\ref{4.22}), we get $c_\alpha=(\sum\limits_{\alpha}c_{\alpha}^2)\bar{c}_\alpha$. Since we assume $\ell'$ is real, then $\sum\limits_{\alpha}c_{\alpha}^2$ is also real. Therefore, we must have $\sum\limits_{\alpha}c_{\alpha}^2=1$ and $c_\alpha$ are real for all $\alpha$. Now we can find a matrix $\textbf{T}\in SO(s)$ s.t. $(c_1,\cdots,c_s)\;\textbf{T}=(1,0,\ldots,0)$. Then the immersion (\ref{4.19})
is congruent to (\ref{4.20}) by choosing the isometry $\textbf{T}\otimes I_{m_1}$ of $\textbf{HP}^{n}$.

If $\sum\limits_{\alpha}c_{\alpha}^2=0$, i.e. $\ell'=0$, then the choosing of $\textbf{z}$ is still allowed a $U(1)$-transformation. Similar to the proof of Proposition \ref{proposition4.1}, the immersion (\ref{4.19}) is congruent to
\begin{equation*}
[a,b]\mapsto[\cos t\;\phi_{1,m_1},\texttt{i}\sin t\;\phi_{1,m_1},0,\ldots,0],
\end{equation*}
where $t\in[0,\pi)$. It follows from $\sum\limits_{\alpha}c_{\alpha}^2=0$ that $t=\pi/4$. Thus we complete the proof. \hfill$\Box$

\vspace{0.3cm}
Next, we calculate the Gauss curvatures of the immersions (\ref{4.20}) and (\ref{4.21}). For the immersion (\ref{4.20}), by (\ref{2.16}), (\ref{4.4}) and straightforward computation, we obtain
$|X|^2+|Y|^2=|\textbf{z}B|^2=a^2_{3,n_1}.$
So, by (\ref{4.6}), the Gauss curvature is
\begin{equation}\label{4.23}
K=\frac{4}{a^2_{3,n_1}}=\frac{16}{(n_1+1)^2-4}.
\end{equation}

For the immersion (\ref{4.21}), by $\ell'=0$, we have $|X|^2+|Y|^2=|\textbf{z}A|^2+|\textbf{z}B|^2=a^2_{1,n_1}+a^2_{3,n_1}.$
Then the Gauss curvature is
\begin{equation}\label{4.24}
K=\frac{4}{a^2_{1,n_1}+a^2_{3,n_1}}=\frac{8}{(n_1+1)^2-2}.
\end{equation}

\vspace{0.3cm}

\textbf{(II.II)}\hspace{0.3cm} If the associate representation $\rho$ takes form of (\ref{4.18}), we split $$\textbf{z}=\sum\limits_{\alpha=1}^{s_1}c_{\alpha}\textbf{u}_{1,n_1}+\sum\limits_{\beta=s_1+1}^sc_{\beta}\textbf{u}_{1,n_2},$$
where $c_\alpha,c_\beta\in\textbf{C}$ and $\sum\limits_{\alpha}c_{\alpha}\bar{c}_\alpha+\sum\limits_{\beta}c_{\beta}\bar{c}_\beta=1$. Then the corresponding immersion is given by
\begin{equation}\label{4.25}
[a,b]\mapsto[c_1\phi_{1,m_1},\ldots,c_{s_1}\phi_{1,m_1},
c_{s_1+1}\phi_{1,m_2},\ldots,c_{s}\phi_{1,m_2}].
\end{equation}

\begin{prop}\label{proposition4.3}
Up to an isometry of $\textbf{\emph{HP}}^{n}$, the minimal immersion \emph{(\ref{4.25})}
is congruent to
\begin{equation}\label{4.26}
f_{m_1,m_2}:[a,b]\mapsto\bigg[\sqrt{\frac{m_1}{m_1+m_2}}\;\phi_{1,m_1},\texttt{\emph{i}}
\sqrt{\frac{m_2}{m_1+m_2}}\;\phi_{1,m_2},
0,\ldots,0\bigg]
\end{equation}
when $m_1+m_2$ is an even,
or
\begin{equation}\label{4.27}
f'_{m_1,m_2}:[a,b]\mapsto\bigg[\sqrt{\frac{m_1}{m_1+m_2}}\;\phi_{1,m_1},
\sqrt{\frac{m_2}{m_1+m_2}}\phi_{1,m_2},0,\ldots,0\bigg]
\end{equation}
when $m_1+m_2$ is an odd, respectively, where $m_1<m_2$.
\end{prop}

\emph{Proof}. For this subcase, the minimality condition (\ref{4.14}) can be written as
\begin{eqnarray}
(a_{1,n_1}^2-p_1)c_\alpha=(-1)^{m_1+1}\ell'a_{1,n_1}\bar{c}_\alpha, \label{4.28} \\
(a_{1,n_2}^2-p_1)c_\beta=(-1)^{m_2+1}\ell'a_{1,n_2}\bar{c}_\beta, \label{4.29}
\end{eqnarray}
for all $\alpha,\ \beta$, where
\begin{equation}\label{4.30}
\ell'=(-1)^{m_1+1}a_{1,n_1}\sum\limits_{\alpha}c_{\alpha}^2+(-1)^{m_2+1}a_{1,n_2}\sum\limits_{\beta}c_{\beta}^2.
\end{equation}
We first claim that $\ell'\neq 0$. Otherwise, if $\ell'= 0$, then from (\ref{4.28}) and (\ref{4.29}), we get $p_1=a_{1,n_1}^2=a_{1,n_2}^2$, i.e. $m_1=m_2$, which leads to a contradiction.
Multiplying by $\bar{c}_{\alpha}$ and $\bar{c}_{\beta}$ on both sides of (\ref{4.28}) and (\ref{4.29}), respectively, and summating for all $\alpha$ and $\beta$, we obtain
\begin{eqnarray}
(a_{1,n_1}^2-p_1)\sum_{\alpha}c_\alpha\bar{c}_{\alpha}=(-1)^{m_1+1}\ell'a_{1,n_1}\sum_{\alpha}\bar{c}^2_\alpha, \label{4.31} \\
(a_{1,n_2}^2-p_1)\sum_{\beta}c_\beta\bar{c}_{\beta}=(-1)^{m_2+1}\ell'a_{1,n_2}\sum_{\beta}\bar{c}^2_\beta, \label{4.32}
\end{eqnarray}
By using the assumption $\ell'$ is real and adding (\ref{4.31}) and (\ref{4.32}), we have
$$p_1=a_{1,n_1}^2\sum_{\alpha}c_\alpha\bar{c}_{\alpha}+a_{1,n_2}^2\sum_{\beta}c_\beta\bar{c}_{\beta}-(\ell')^2,$$
so $p_1$ is also real. Therefore, it follows from (\ref{4.28}) and (\ref{4.29}) that $c_\alpha=\pm\bar{c}_\alpha$ for all $\alpha$ and $c_\beta=\pm\bar{c}_\beta$ for all $\beta$. Moreover, we have
$$a_{1,n_1}^2-p_1=\pm(-1)^{m_1+1}\ell'a_{1,n_1},\ \ a_{1,n_2}^2-p_1=\pm(-1)^{m_2+1}\ell'a_{1,n_2}.$$
Without lose of generality, we first consider $c_\alpha=\bar{c}_\alpha$ for all $\alpha$ and $c_\beta=\bar{c}_\beta$ for all $\beta$, i.e. $c_\alpha,c_\beta$ are all real. Then $\sum\limits_{\alpha}c_{\alpha}^2+\sum\limits_{\beta}c_{\beta}^2=1$ and
\begin{eqnarray*}
a_{1,n_1}^2-p_1=(-1)^{m_1+1}\ell'a_{1,n_1}, \ \ a_{1,n_2}^2-p_1=(-1)^{m_2+1}\ell'a_{1,n_2},
\end{eqnarray*}
which gives
$$a_{1,n_1}^2- a_{1,n_2}^2=\ell'\big((-1)^{m_1+1}a_{1,n_1}-(-1)^{m_2+1}a_{1,n_2})\big),$$
so
\begin{eqnarray}\label{4.33}
\ell'=(-1)^{m_1+1}a_{1,n_1}+(-1)^{m_2+1}a_{1,n_2}.
\end{eqnarray}
Combining (\ref{4.30}) and (\ref{4.33}), we obtain
$$\sum_{\beta}c_\beta^2=(-1)^{m_1+m_2+1}\frac{a_{1,n_2}}{a_{1,n_1}}\sum_{\alpha}c_\alpha^2>0,$$
which follows that $m_1+m_2$ is an odd number and
\begin{eqnarray}
\sum_{\alpha}c_\alpha^2=\frac{a_{1,n_1}}{a_{1,n_1}+a_{1,n_2}}=\frac{m_1}{m_1+m_2},\ \ \sum_{\beta}c_\beta^2=\frac{a_{1,n_2}}{a_{1,n_1}+a_{1,n_2}}=\frac{m_2}{m_1+m_2}.\label{4.34}
\end{eqnarray}
Therefore, by similar argument as the proof of Proposition \ref{proposition4.1}, the immersion (\ref{4.25}) is congruent to (\ref{4.27}).

Next, we only need to consider $c_\alpha=\bar{c}_\alpha$ for all $\alpha$ and $c_\beta=-\bar{c}_\beta$ for all $\beta$, i.e. $c_\alpha$ are all real and $c_\beta$ are all purely imaginary. Then $\sum\limits_{\alpha}c_{\alpha}^2-\sum\limits_{\beta}c_{\beta}^2=1$ and
\begin{eqnarray*}
a_{1,n_1}^2-p_1=(-1)^{m_1+1}\ell'a_{1,n_1}, \ \ a_{1,n_2}^2-p_1=-(-1)^{m_2+1}\ell'a_{1,n_2},
\end{eqnarray*}
which gives
$$a_{1,n_1}^2- a_{1,n_2}^2=\ell'\big((-1)^{m_1+1}a_{1,n_1}+(-1)^{m_2+1}a_{1,n_2})\big),$$
so
\begin{eqnarray}\label{4.35}
\ell'=(-1)^{m_1+1}a_{1,n_1}-(-1)^{m_2+1}a_{1,n_2}.
\end{eqnarray}
Combining (\ref{4.30}) and (\ref{4.35}), we obtain
$$\sum_{\beta}c_\beta^2=(-1)^{m_1+m_2+1}\frac{a_{1,n_2}}{a_{1,n_1}}\sum_{\alpha}c_\alpha^2<0,$$
which follows that $m_1+m_2$ is an even number and
\begin{eqnarray}\label{4.36}
\sum_{\alpha}c_\alpha^2=\frac{a_{1,n_1}}{a_{1,n_1}+a_{1,n_2}}=\frac{m_1}{m_1+m_2},\ \ \sum_{\beta}c_\beta^2=-\frac{a_{1,n_2}}{a_{1,n_1}+a_{1,n_2}}=-\frac{m_2}{m_1+m_2}.
\end{eqnarray}
Therefore, by similar argument as the proof of Proposition \ref{proposition4.1}, the immersion (\ref{4.25}) is congruent to (\ref{4.26}). Thus we complete the proof. \hfill$\Box$

\vspace{0.2cm}

Now, we calculate the Gauss curvatures of the immersions (\ref{4.26}) and (\ref{4.27}). For the immersion (\ref{4.27}), by (\ref{2.16}), (\ref{4.4}), (\ref{4.34}) and straightforward computation, we obtain
\begin{eqnarray*}
|X|^2=|\textbf{z}A-\ell\textbf{z}|^2&=&a^2_{1,n_2}\sum_{\alpha}c_{\alpha}^2+a^2_{1,n_1}\sum_{\beta}c_{\beta}^2\\
                                    &=&m_2^2\frac{m_1}{m_1+m_2}+m_1^2\frac{m_2}{m_1+m_2}=m_1m_2,
\end{eqnarray*}
and
\begin{eqnarray*}
|Y|^2=|\textbf{z}B|^2&=&a^2_{3,n_1}\sum_{\alpha}c_{\alpha}^2+a^2_{3,n_2}\sum_{\beta}c_{\beta}^2\\
                                    &=&(m_1-1)^2\frac{m_1}{m_1+m_2}+(m_2-1)^2\frac{m_2}{m_1+m_2}=m_1^2-m_1m_2+m_2^2-1.
\end{eqnarray*}
So by (\ref{4.6}), the Gauss curvature is
\begin{equation}\label{4.37}
K=\frac{4}{m_1^2+m_2^2-1}.
\end{equation}
For the immersion (\ref{4.26}), we can obtain the expression of the Gauss curvature by similar calculation, which is also given by (\ref{4.37}).

\vspace{0.3cm}

In summary, we completely classify the minimal homogeneous 2-spheres in quaternionic projective space $\textbf{HP}^n$. That is

\begin{theo}\label{theorem4.4}
Let $f:S^2\longrightarrow\textbf{\emph{HP}}^n$ be a linearly full homogeneous minimal immersion. Then, in terms of homogeneous coordinates, $f$ is congruent to one of the following\emph{:}

\emph{(1)} $f_\lambda:[a,b]\mapsto [\phi_{\lambda,n+1}],$ $\lambda\in\{3,5,\ldots,2n+1\}$, and the Gauss curvature $K=8/[(2n+2)^2-(\lambda^2+1)];$

\emph{(2)} $f_1:[a,b]\mapsto[\phi_{1,n+1}]$ and the Gauss curvature $K=4/[n(n+2)];$

\emph{(3)} $f_{\lambda,m,t}:[a,b]\mapsto[\cos t\phi_{\lambda,m},\texttt{\emph{i}}\sin t\phi_{\lambda,m}]$ for some positive weight $\lambda\in\{3,5,\ldots,n\}$, $t\in(0,\pi/2)$, $2m=n+1$, and the Gauss curvature $K=8/[(n+1)^2-(\lambda^2+1)];$

\emph{(4)} $f_{m_1,m_2}:[a,b]\mapsto[\sqrt{m_1/(m_1+m_2)}\;\phi_{1,m_1},
\texttt{\emph{i}}\sqrt{m_2/(m_1+m_2)}\;\phi_{1,m_2}]$ for some positive $m_1\leq m_2$ so that $m_1+m_2=n+1$ is an even, and the Gauss curvature $K=4/(m_1^2+m_2^2-1);$

\emph{(5)} $f'_{m_1,m_2}:[a,b]\mapsto[\sqrt{m_1/(m_1+m_2)}\;\phi_{1,m_1},
\sqrt{m_2/(m_1+m_2)}\;\phi_{1,m_2}]$ for some positive $m_1<m_2$ so that $m_1+m_2=n+1$ is an odd, and the Gauss curvature $K=4/(m_1^2+m_2^2-1).$

\end{theo}
\emph{Proof}. From the Proposition \ref{proposition4.1}, \ref{proposition4.2} and \ref{proposition4.3}, we know that $f$ is congruent to $f_\lambda$, or $f_1$, or $f_{\lambda,m,t}$, or $f_{m_1,m_2}$, or $f'_{m_1,m_2}$. The expressions for the Gauss curvature $K$ follows from (\ref{4.12}), (\ref{4.23}), (\ref{4.24}) and (\ref{4.37}), respectively. \hfill$\Box$

\vspace{0.3cm}

We should point out the congruence of $f_{\lambda,m,t}$ w.r.t. the parameter $t$, that is
\begin{prop}\label{proposition4.5}
$f_{\lambda,m,t_1}$ is not congruent to $f_{\lambda,m,t_2}$ for any $t_1,t_2\in(0,\pi/2)$, $t_1<t_2$.
\end{prop}
\emph{Proof}. Suppose that $f_{\lambda,m,t_1}$ is congruent to $f_{\lambda,m,t_2}$ for $t_1,t_2\in(0,\pi/2)$, then there exists
fixed $p\in Sp(1)$, $P\in Sp(2m)$ such that
\begin{equation}\label{4.38}
p(\cos t_1,\texttt{i}\sin t_1)\rho_m\oplus\rho_m(g)=(\cos t_2,\texttt{i}\sin t_2)\rho_m\oplus\rho_m(g)P
\end{equation}
for all $g\in SU(2)$. Set $p=p'+p''\texttt{j}$ and $P=(p'_{AB})+(p''_{AB})\texttt{j}$. Since $\{\Lambda_{AB}\;|\;0\leq A,B\leq 2m-1\}$ are linearly independent in $L^2(SU(2))$ by Peter-Weyl's theorem,
we obtain
\begin{equation}\label{4.39}
p''=p''_{AB}=p'_{s\;l}=p'_{s\;m+l}=p'_{m+s\;l}=p'_{m+s\;m+l}=0,\;\;s\neq l,
\end{equation}
and
\begin{equation}\label{4.40}
(p'+\bar p')\cos t_1=2\cos t_2\; p'_{l\;l},\;\;\;
(p'-\bar p')\cos t_1=2\texttt{i}\sin t_2\; p'_{m+l\;l},
\end{equation}
\begin{equation}\label{4.41}
(p'-\bar p')\sin t_1=-2\texttt{i}\cos t_2\; p'_{l\;m+l},\;\;\;
(p'+\bar p')\sin t_1=2\sin t_2\; p'_{m+l\;m+l},
\end{equation}
by comparing the coefficients of $\Lambda_{AB}$. The identities (\ref{4.40}) and (\ref{4.41}) imply that $t_1=t_2$. \hfill$\Box$

\vspace{0.3cm}

Base on the study of homogeneous harmonic maps into projective space, Ohnita \cite{[Oh]} conjectured that all proper minimal constant curved minimal two-spheres in $\textbf{HP}^n$ are $\{f_\lambda:\lambda=1,3,\ldots,2n+1\}$, up to a rigidity of $\textbf{HP}^n$.
Notice that $\textbf{HP}^n$ is a rank one compact symmetric space,
together with Theorem \ref{theorem4.4},
we propose the following conjecture, still called Ohnita's conjecture:
\begin{conj}\label{conjecture4.5}
Let $f$ be a proper minimal constant curved immersion from $S^2$ into $\textbf{\emph{HP}}^n$, then $f$ is congruent to one of $f_\lambda$, $f_{\lambda,m,t}$, $f_{m_1,m_2}$, $f'_{m_1,m_2}$, where $\lambda\in\{1,3,\ldots,2n+1\}$.
\end{conj}

\vspace{0.3cm}
\noindent\textbf{Acknowledgments}. This work was supported by NSFC (Nos. 11471299, 11401481, 11331002).

\noindent Jie Fei,\\
Department of Mathematical Sciences, Xi'an Jiaotong-Liverpool University,
Suzhou 215123, P. R. China. E-mail: jie.fei@xjtlu.edu.cn

\vspace{0.3cm}
\noindent Chiakuei Peng,\\
School of Mathematical Sciences, University of Chinese Academy of
Sciences, Beijing 100049, P. R. China. Email: pengck@ucas.ac.cn.

\vspace{0.3cm}
\noindent Xiaowei Xu, Corresponding author\\
 Department of Mathematics,  University of Science and Technology of China, Hefei, 230026, Anhui
province, P. R. China;\\ and Wu Wen-Tsun Key Laboratory of Mathematics, USTC, Chinese Academy of
Sciences, Hefei, 230026, Anhui, P. R. China. Email: xwxu09@ustc.edu.cn\\


\begin{thebibliography}{2}



\bibitem{[BO]} S. Bando and Y. Ohnita, {\it Minimal $2$-spheres with constant curvature in $\mathbb{C}P^n$}, J.Math.Soc.Japan. 3(1987), 477-487.

\bibitem{[BD]} T. Br$\ddot{o}$cker and T. tom Dieck, {\it Representations of compact Lie groups}, GTM98, Springer-Verlag, New York, 1985.


\bibitem{[BJRW]} J. Bolton, G.R. Jensen, M. Rigoli and L. M.
Woodward, {\it On conformal minimal immersions of $S^2$ into
$\mathbb{C }P^n$}, Math. Ann. 279(1988), 599-620.



\bibitem{[CW2]} S.S. Chern and J.G. Wolfson, {\it Harmonic maps of the two-sphere into
a complex Grassmann manifold II}, Ann. of Math. 125 (1987), 301-335.

\bibitem{[Ca]} E. Calabi, {\it Minimal immersions of surfaces in Euclidean
spheres}, J. Diff. Geom. 1(1967), 111-126.



\bibitem{[DW]} M.P.do Carmo, N.R.Wallach, {\it Minimal immersion of
spheres into spheres}, Ann. of Math. 93(1971), 43-62.



\bibitem{[DHZ]} L.Delisle, V.Hussin and W.J.Zakrzewski, {\it Constant curvature solutions of Grassmannian sigma
models: \emph{(1)} Holomorphic solutions}, Journal of Geometry and Physics 66(2013), 24-33.

\bibitem{[DHZ2]} L.Delisle, V.Hussin and W.J.Zakrzewski, {\it Constant curvature solutions of Grassmannian sigma
models: \emph{(2)} Non-holomorphic solutions}, Journal of Geometry and Physics 71(2013), 1-10.



\bibitem{[FJXX]} J.Fei, X.X.Jiao, L.Xiao, X.W. Xu, {\it On the Classification of
homogeneous 2-spheres in complex Grassmannians}, Osaka
J. Math., 50(2013), 135-152.

\bibitem{[HJ]} L.He and X.X.Jiao {\it Classification of conformal minimal immersions of constant curvature from $S^2$ to $HP^2$}, Math. Ann., 359(2014), No.3-4, 663-694.

\bibitem{[LY]} Zh.Q. Li and Zh.H. Yu, {\it Contant curved minimal
2-spheres in $G(2,4)$}, Manuscripta Math. 100(1999), 305-316.



\bibitem{[JP]} X.X.Jiao and J.G.Peng, {\it Classification of
holomorphic two-spheres with constant curvature in the complex
Grassmannians $G_{2,5}$}, Diff. Geom. Appl., 20(2004), 267-277.

\bibitem{[Oh]} Y.Ohnita, {\it Homogeneous harmonic maps into projective space}, Tokyo J.Math. 13(1)1990, 87-116.

\bibitem{[PX]} C.K.Peng and X.W.Xu, {\it Classification of minimal homogeneous two-spheres in the complex Grassmann manifold G(2,n)},  J. Math. Pures Appl., (9)103(2015), No.2, 374-399.

\bibitem{[PJX]} C.K.Peng, J.Wang and X.W.Xu, {\it Minimal two-spheres with constant curvature in the complex hyperquadric}, J. Math. Pures Appl., (9)106(2016), No.3, 453-476.





\end{thebibliography}
\end{document}